\documentstyle[12pt,amsfonts]{article}

\newcommand{\Pic}{{\rm Pic}}

\newcommand{\Spec}{{\rm Spec \,}}

\newcommand{\cD}{{\cal D}}

\newcommand{\cJ}{{\cal J}}

\newcommand{\cO}{{\cal O}}

\newcommand{\A}{{\Bbb A}}

\newcommand{\C}{{\Bbb C}}

\newcommand{\E}{{\Bbb E}}

\newcommand{\G}{{\Bbb G}}

\newcommand{\N}{{\Bbb N}}

\renewcommand{\P}{{\Bbb P}\hspace{.05em}}

\newcommand{\Z}{{\Bbb Z}}

\newcommand{\iC}{{\rm C\hspace{-.4em}l\hspace{.2261em}}}

\title{\bf $\C^*$-extensions of tori, higher Chow groups and applications
 to incidence equivalence relations for algebraic cycles }

\author{by Stefan {\sc M\"uller-Stach}, Essen}

\date{January 1995}

\begin{document}

\maketitle

\ \\ \\ \\
{\bf Abstract:}{\small Let $X$ be a smooth projective variety over an
algebraically closed field $k$. We repeat Bloch's construction of a
$\G_m$-biextension (torseur) $\E$ over
$CH^{p}_{{\rm hom}} (X) \times CH^{q}_{{\rm hom}} (X)$ for $p+q=\dim(X)+1$.
First we show that in characteristic zero $\E$ comes via pullback from the
Poincar\'e biextension $\P$ over the corresponding product of intermediate
Jacobians - as conjectured by Bloch and Murre. Then the relations
between $\E$ and various equivalence relations for algebraic cycles are
studied. In particular we reprove Murre's theorem stating that Griffiths'
conjecture holds for codimension $2$ cycles, i.e. every $2$-codimensional
cycle which is algebraically and incidence equivalent to zero has torsion
Abel-Jacobi invariant.} 
\newpage
\noindent
\underline{{\bf \S \ 0. Introduction}} \\
\ \\ \\
Let $X$ be a compact K\"ahler manifold of dimension $d$. If we fix
two integers $0 \leq p, q \leq d$ with $p + q = d +1 $, then the two
Griffiths intermediate Jacobians $\cJ^p
(X) $ and $\cJ^q (X)$ are mutually dual. The
Poincar\'{e} line bundle $P$ on $\cJ^p (X) \times
\cJ^q (X)$ minus its zero section, denoted by $\P$, is a
$\C^*$-torseur and in fact a biextension in the sense of
Grothendieck (SGA 7) and [Mumford 69]. \\
\ \\
Any two disjoint homologically trivial cycles $Z \in
CH^{p}_{{\rm hom}} (X)$, $W \in CH^{q}_{{\rm hom}} (X)$ define
an element $<Z, W> \in \P$, projecting onto the element
$(\psi_p (Z)$, $\psi_q (W)) \in \cJ^p (X) \times \cJ^q
(X)$, where $\psi$ is the Abel-Jacobi homomorphism. \\
\ \\
Bloch has given another construction for a $\G_m$-torseur
$\E$ over $CH^{p}_{{\rm hom}} (X) \times CH^{q}_{{\rm hom}} (X)$
for a projective manifold $X$ over an arbitrary field [Bloch
89]. He also conjectures, that $\E$ may be obtained as the
pullback of $\P$ via the Abel-Jacobi map in characteristic zero.
We prove this:\\
\ \\
\proclaim Theorem 1. Let $X$ be a smooth projective variety over $\C$. Then
$\E$ is the pullback of $ \P $ as torseurs
under the Abel-Jacobi homomorphism.
\par
\ \\
The idea of proof is that the essential point in Bloch's
construction factors through Deligne cohomology. \\
\ \\
Theorem 1 can be applied to shed some light on the relation
between Abel-Jacobi equivalence and incidence equivalence of
algebraic cycles. A cycle $W \in CH^{q}_{{\rm hom}} (X)$ is
called {\it incidence equivalent} to zero $(W \sim_{{\rm inc}} \ 0)$,
if for all pairs $(T,B)$, with
a smooth projective variety $T$ and a cycle $B \in
CH^{d+1-q} (T \times X)$, the divisor
$B(W) = (pr_T)_* (B \cap pr^{*}_{X} (W))$ is linearly equivalent
to zero on $T$ (see [Griffiths 69]). \\
\ \\
It is easy to show that $\psi (W) = 0$ implies $W \sim_{{\rm inc}}
\ 0$. Griffiths conjectured, that the opposite is true modulo
torsion if $W$ is assumed to be algebraicly equivalent to zero.
It was known that this follows from the generalized Hodge
conjecture of Grothendieck. Murre ([Murre 85]) has proved
the conjecture for codimension two cycles. We show here, that theorem 1 has
something to do with this problem. Denote by $\E_W^{\rm alg}$ the
restriction of $\E$ to the fiber over $W \in CH^{q}_{{\rm alg}}
(X)$:
$$
1 \to \C^* \to \E_W^{\rm alg} \to CH^{p}_{{\rm alg}} (X) \to 0
$$
We show: \\
\ \\
\proclaim Theorem 2.
Let $X$ be smooth projective over $\C$ and $W$ be algebraically equivalent
to zero. Then
$\E_W^{\rm alg}$ is a split extension if and only if $W \sim_{{\rm inc}} \ 0$.
\par
\ \\
Here the idea is that incidence equivalence reduces somehow to the situation
of points and divisors, which is better understood than the general case.
This result brings together the a priori different definitions of
incidence equivalence occuring in the literature. Theorem 2
implies a condition for Griffiths' conjecture: \\
\ \\
\proclaim Theorem 3.
Let $X$ be as above and
$W \in CH^{q}_{{\rm alg}} (X)$ be such that $\psi (N \cdot W)$ is
contained in the dual of $\cJ^{p}_{{\rm alg}} (X) = \psi
(CH^{p}_{{\rm alg}} (X))$ (an abelian subvariety) for some $N \in \N$.
Then $W \sim_{{\rm inc}} \ 0$ implies that
$\psi (W)$ is torsion in $\cJ^q (X)$.
\par
\ \\
The condition follows again from Grothendieck's generalized
Hodge conjecture. But it holds always for codimension two cycles if and
therefore we get another proof of Murre's theorem ([Murre85]): \\
\ \\
\proclaim Corollary(Murre's theorem). \\
Griffiths'conjecture holds for codimension two cycles on a smooth projective
manifold $X$ over $\C$, i.e. for every cycle $W$ on $X$, algebraicly and
incidence equivalent to zero, $\psi(W)$ is torsion in $\cJ^2(X)$.
\par
\ \\
It is easy to see that the proof of Griffiths' conjecture
would go through in any codimension $p$, if the adjoint $\Lambda$ of the
Lefschetz operator were again algebraic. This can be verified
for several types of varieties e.g. rational-like varieties and abelian
varieties and complete intersections in these. For a special case see
[Griffiths-Schmid75].\\
The paper is organized as follows: \\
\ \\
\S \ 1. The Poincar\'{e} biextension. \\
\S \ 2. Facts about higher Chow groups. \\
\S \ 3. Bloch's construction. \\
\S \ 4. The pullback theorem. \\
\S \ 5. Abel-Jacobi versus incidence equivalence. \\
\ \\
It is a pleasure to thank Jacob Murre for mentioning this problem
to me and showing me his unpublished notes ([Murre])
during a visit at Leiden. Furthermore I am grateful to A. Collino, H. Esnault
and D. Huybrechts for many helpful remarks and improvements. \\
\newpage
\noindent
\underline{{\bf \S \ 1. The Poincar\'{e} biextension}} \\
\ \\ \\
Let us first discuss briefly the notion of biextension. We refer to
[Mumford 69] and SGA 7 for further details: Let $A,B,C$ be abelian groups.
A {\sl biextension} of $B \times C$ by $A$ is a set $E$ on which $A$ acts
freely, together with a quotient map $\pi: E \to B \times C$ and two laws
of composition $+_1:E \times_B E \to E$, $+_2:E \times_C E \to E$ subject to
the following conditions:\\
(1) For all $b \in B$, $\pi^{-1}(b \times C)$ is an abelian group under
$+_1$ and the following sequence is exact:
$$ 0 \to A \to \pi^{-1}(b \times C) \to C \to 0 $$
\\
(2) For all $c \in C$, $\pi^{-1}(B \times c)$ is an abelian group under
$+_2$ and the following sequence is exact:
$$ 0 \to A \to \pi^{-1}(B \times c) \to B \to 0 $$
\\
(3) $+_1,+_2$ are compatible, i.e. for suitable $w,x,y,z \in E$ we have
$$(w+_1 x)+_2 (y +_1 z) = (w +_2 y) +_1 (x+_2 z) $$ \\
There is one prominent example which caused this definition to exist, namely
the {\it Poincar\'e biextension}:
\ \\
Let $T$ be a compact complex torus of dimension $n$. There is a
natural line bundle $P$ (the Poincar\'{e} line bundle) on $T
\times T^\vee$, where $T^\vee = \Pic^0(T)$. As an element of $\Pic
(T \times T^\vee)$ it is usually normalized by the two conditions \\
\ \\
(i) $P |_{\{ 0\} \times T^\vee} \cong \cO_{T^\vee}$ \\
(ii) $P |_{T \times \{ \lambda \}}$ represents $\lambda \in
\Pic^0 (T) = T^\vee$ \\
\ \\
$P$ is unique under these assumptions. Let $\P = P
\setminus$ zero section. Clearly $\C^\ast$ acts freely on $\P$.
$\P$ is called {\it Poincar\'{e}-biextension} associated to $T$.
The projection map
$p: \P \to T \times T^\vee$ makes $\P$ a $\C^*$-torseur over $T
\times T^\vee$. For every $\lambda \in T^\vee$, the inverse image
$p^{-1} (T \times \{ \lambda \})$ is denoted by $\P_{\lambda}$
and sits in the exact sequence
$$
0 \to \C^* \to \P_{\lambda} \to T \to 0.
$$
Then $\P_{\lambda}$ is an extension of abelian groups. It is
well known that
$$T^\vee = \Pic^0 (T) \cong {\rm Ext}^1 (T, \C^*)$$
in the category of complex analytic
groups. $\P_{\lambda}$ is exactly the extension of $T$ by $\C^*$
given by $\lambda$ in this isomorphism. If $D$ is any divisor on $T$ with
$c_1 (\cO (D)) = \lambda$ then we will write also $\P_D$ instead
of $\P_{\lambda}$ and note that $\P_D$ depends only on the class
on $D$ in $\Pic^0 (T)$. Now let us consider a special case: Let
$X$ be a K\"ahler manifold of dimension $d$ and $p$ some integer.
If we let $T = \cJ^p (X)$ the $p$-th intermediate Jacobian of $X$
$$
\cJ^p (X) = \frac{H^{2 p-1} (X, \C)}{F^p \oplus H^{2p-1} (X,\Z)}
$$
then it follows by using Poincar\'{e} duality, that $T^\vee$ is given
by $\cJ^q (X)$ where $q = d +1 -p$. Let
$$
\P \to T \times T^\vee = \cJ^p (X) \times \cJ^q (X)
$$
be the Poincar\'{e}-biextension as defined above. If $W \in
CH^{q}_{{\rm hom}} (X)$ with $\lambda =\psi(W) \in \cJ^q (X)$ and
$\P_W : = \P_{\lambda}$ we remark that $\P_W$
depends only on the Abel-Jacobi equivalence class of $W$,
i.e. $W$ is Abel-Jacobi equivalent
to zero if and only if $\P_W$ splits as an extension in ${\rm
Ext}^1 (T, \C^* )$. \\
\ \\

\ \\ \\ \\
\underline{{\bf \S \ 2. Facts about higher Chow groups}} \\
\ \\ \\
Literature: [Bloch 86], [Murre] for this chapter. Before we
repeat Bloch's construction, we recall some properties of {\it
higher Chow groups} as defined by Spencer Bloch. \\
\ \\
Let $X$ be a quasi-projective variety over $\C$. Denote by $\Delta^n
= \Spec (k [T_0 , \ldots ,T_n ]/<{\Sigma T_i =1}> )$ which is
isomorphic to affine space $\A^{n}_{\iC}$.
There are $n+1$ natural faces isomorphic
to $\Delta^{n-1}$ contained in $\Delta^n$, defined by the vanishing of one
of the coordinate functions $t_i$. Let $Z^r (X, n)$ be the free
abelian group of cycles $Z \subset X \times \Delta^n$ of codimension $r$
meeting all faces $X \times \Delta^m$ $(m < n)$ properly. Then
$CH^r(X,n)$ is defined as the $n$-th homology group of the complex
$$
Z^r (X;\cdot) = ( \ldots Z^r (X, n+1) {\buildrel \partial \over \to}
Z^r (X,n)
{\buildrel \partial \over \to} Z^r (X, n-1) \to \ldots \to Z^r (X, 0)).
$$
The maps $\partial$ are given by alternating sums of restriction
maps to faces. We will need the following properties: \\
\ \\
(1) $CH^* (X, *)$ are covariant (contravariant)
functorial for proper (flat) morphisms. \\
\ \\
(2) If $W \subset X$ is closed, we have a long exact sequence
$$
\ldots CH^* (X, W, n) \to CH^* (X,n) \to CH^* (W, n) \to CH^*
(X, W, n-1) \ldots
$$
(3) $CH^* (X,0) = CH^* (X)$ (ordinary Chowgroups).\\
\ \\
(4) There is a product for $X$ smooth
$$
CH^p (X, q) \otimes CH^r (X, s) \to CH^{p+r} (X, q+s).
$$
(5) There are natural maps to Deligne-Beilinson cohomology, defined in Bloch's
article in Contemporary Math. {\bf 58} (1986):
$$
c_\cD : CH^p (X, q) \to H^{2p-q}_{\cD} (X, \Z (p)).
$$
(6) If $X$ is proper (not necessarily smooth) and $\dim X
= d$, then there is a natural {\it surjective} homomorphism
$$
\epsilon : CH^{d+1} (X, 1)  \to \C^*
$$
factoring through Deligne-Beilinson cohomology
$$
\epsilon : CH^{d+1} (X, 1) {\buildrel c_\cD \over \to}
H^{2d+1}_{\cD} (X, \Z (d+1)) \to \C^*
$$
obtained as follows: Let $\pi: X \to \Spec(\C)$ be the natural morphism. Then
$\epsilon$ is given by the the direct image map
$\pi_\ast: CH^{d+1}(X,1) \to CH^1(\Spec(\C),1)$, since it is a
straightforward exercise to show that
$$CH^1(\Spec(\C),1) \cong H^1_{\cD}(\Spec(\C),\Z(1)) \cong \C/\Z(1)$$
via the classes in (5) and the latter group can be identified with
$\C^*$. We claim:\\
(1) $\epsilon$ is surjective.\\
(2) It factors through Deligne-Beilinson cohomology.\\
\underline{Proof for (1) and (2)}\\
(1): By definition $CH^1(\Spec(\C),1)=Z^1(\Spec(\C),1)/{\rm Im}(\partial)$,
since $Z^1(\Spec(\C),0)=0$. In the same way
$CH^{d+1}(X,1)=Z^{d+1}(X,1)/{\rm Im}(\partial)$.
Hence in both groups the elements
are represented by finite sums of points. It is clear that the
map $  Z^{d+1}(X,1) \to Z^1(\Spec(\C),1)$ induced by $\pi_\ast$ is
surjective. Hence we deduce (1).\\
(2): By the functorial properties of the Deligne-Beilinson classes we get
a commutative diagram:
$$\matrix{ CH^{d+1}(X,1) & {\buildrel \pi_\ast \over \to} &
CH^1(\Spec(\C),1) \cr
\downarrow c_{\cD} & & \downarrow c_{\cD}  \cr
H^{2d+1}_{\cD}(X,\Z(d+1)) & {\buildrel \pi_\ast \over \to} &
H^1_{\cD}(\Spec(\C),\Z(1))     }
$$
The right vertical arrow being the identity we get (2). \hfill $\square$ \\
\newpage
\noindent
\underline{{\bf \S \ 3. Bloch's Construction}} \\
\ \\ \\
For this chapter we refer to [Bloch 89].\\
Let us describe Bloch's construction in detail: Let $X$ be smooth and
projective over $\C$ with $\dim X =d$, and $W \in Z_{m-1} (X) =
Z^{d+1-m} (X)$ a cycle of dimension $m-1$, which is homologous
to zero, denoted by $W \sim_{{\rm hom}} \ 0$ as usual. We will
construct an extension
$$
1 \to \C^* \to \E_W \to CH^{m}_{{\rm hom}} (X) \to 0
$$
such that they can be glued together to give
$$
\E = \bigcup \E_W \to CH^{m}_{{\rm hom}} (X) \times
CH^{d+1-m}_{{\rm hom}} (X)
$$
a biextension in the sense of \S \ 1. To do this, Bloch
first considers the map
$$
\Theta_W = \epsilon \circ (\cap W) : CH^m (X, 1) {\buildrel \cap W \over \to}
CH^{d+1} (X, 1) {\buildrel \epsilon \over \to} \C^*
$$
using properties (4) and (6) of higher Chow groups. \\
\ \\
\proclaim Lemma 1.
If $W \sim_{{\rm hom}} \ 0 $, then $\Theta_W \equiv 1$.
\par
\ \\
{\bf Proof.} Recall the following fact from Deligne-Beilinson cohomology
([Esnault-Viehweg 88]): \\
$\cJ^* \subset H^{*}_{\cD} (X, \Z (*))$ is a square zero ideal,
where
$$\cJ^* = {\rm Ker} (H^{*}_{\cD} (X, \Z (*)) \to
H^{*}_{\rm Betti} (X, \Z(*)))$$
Now look at the diagram
$$
\matrix{
CH^m (X, 1) &{\buildrel {\cap W} \over \to}& CH^{d+1} (X, 1)
& {\buildrel \epsilon \over \to} & \C^* \cr
\downarrow c_\cD &&\downarrow  c_\cD && || \cr
H^{2m-1}_{\cD} (X, \Z (m)) &{\buildrel {\cap c_\cD(W)}\over \to }&
H^{2d+1}_{\cD} (X,\Z (d+1))& {\buildrel  \epsilon \over \to }& \C^*}
$$
But $W \sim_{{\rm hom}} \ 0$ implies $c_\cD (W) \in \cJ^* (X)$, and
also $H^{2m-1}_{\cD} (X, \Z (m)) \subset \cJ^* (X)$. Hence the
lemma follows. \hfill $\square$ \\
\ \\
To construct $\E_W$ consider the long exact sequence for the
pair $(X, |W|)$:
$$
\ldots \to CH^m (X, 1) \to CH^m (|W|,1) \to CH^m (X, |W|, 0)
\to CH^m (X, 0) \to 0
$$
where $|W| = {\rm supp} (W) $ as usual. Since $|W|$ is proper of
dimension $m-1$, we have a surjective map
$\epsilon : CH^m (|W|,1) \to  \C^*$ by property (6).
But $\Theta_W \equiv 1$, hence $\epsilon$ factors
through
$$
A: = \frac{CH^m (|W|, 1)}{{\rm Im} \  CH^m (X, 1)} \to  \C^* .
$$
Consider the commutative diagram
$$
\matrix{
0 &\to & A &{\buildrel i \over \to} & CH^m (X, |W|, 0) &\to &
CH^m (X) &\to & 0 \cr
&& \downarrow \epsilon && \downarrow && || && \cr
1 &\to& \C^* &\to& \E_W &\to&  CH^m (X) &\to 0}
$$
where we have defined $\E_W$ as
$$
CH^m (X, |W|, 0) / i ({\rm Ker} (\epsilon) )
$$
as the pushout of the upper line via $\epsilon$. By abuse of
language, we also denote by $\E_W$ the restriction
$$
1 \to \C^* \to \E_W \to CH^{m}_{{\rm hom}} (X) \to 0
$$
to $CH^{m}_{{\rm hom}} (X)$.
Bloch proves several results for $\E_W$: \\
If $W$ is rationally equivalent to $W'$, then $\E_W$ is
canonically isomorphic to $\E_{W'}$. In particular, if $W
\sim_{{\rm rat}} \ 0$, then $\E_W$ splits as an extension of
groups. We prove a little bit more in \S \ 4. \\
\ \\
The $\E_W$ may be glued together, to obtain a $\C^*$-torseur (biextension)
$$
\E \to CH^{m}_{{\rm hom}} (X) \times CH^{d+1-m}_{{\rm hom}} (X)
$$
satisfying the axioms of \S \ 1. \\
We refer the reader to [Bloch 89] for this construction, since
we don't use the full torseur $\E$ for our applications. In fact
by our theorem 1 this follows from the statement for $\P$. \\
\ \\ \\ \\
\underline{{\bf \S \ 4. The Pullback Theorem}} \\
\ \\ \\
Here we prove the main result of this paper.
Let $X$ be a smooth projective variety over $\C$. \\
\proclaim Lemma 2.
If $W \sim_{A\cJ} \ 0$ (Abel-Jacobi equivalent to zero), then
$\E_W$ splits.
\par
\ \\
{\bf Proof.} We use the functorial Deligne-Beilinson classes $c_\cD$ from (5),
\S 2. Consider the commutative diagram
$$
\matrix{
A &\to & CH^{m}_{{\rm hom}} (X, |W|, 0) &\to&
CH^{m}_{{\rm hom}} (X) &\to &  0 \cr
\downarrow c_\cD && \downarrow && \downarrow & {\rm Abel-Jacobi} & \cr
H^{2m -1}_{\cD} (|W|; \Z (m)) & \to & H^{2m}_{\cD} (X, |W|; \Z
(m))_{{\rm hom}}&  \to &  \cJ^m (X) &\to & 0}
$$
$\epsilon : A \to  \C^* $ is surjective and by \S 2.(6) factors through
$$
H^{2m-1}_{\cD} (|W|; \Z (m)) \to \bigoplus \C^* {\buildrel \Sigma \over
\to }  \C^*
$$
by summation over all irreducible components of $|W|$. Define
$K:= {\rm Ker} (A \to  \C^* )$ and recalling that $\E_W =
CH^{m}_{{\rm hom}} (X, |W|, 0)/i(K)$ we obtain the
diagram
$$
\matrix{
 &&&&T && T'&& \cr
&&&&  \downarrow & & \downarrow && \cr
0& \to &A/K &\to &\E_W &\to & CH^{m}_{{\rm hom}} (X) &\to &0 \cr
&&\downarrow \epsilon && \downarrow & & \downarrow &{\rm Abel-Jacobi}
& \cr
1& \to& \C^*& \to& \P_W& \to& \cJ^m (X)& \to & 0}
$$
where $\P_W$ is as in \S \ 1, $T = {\rm Ker} (\E_W \to \P_W )$ and
$T' = {\rm Ker} (CH^{m}_{{\rm hom}} (X) \to \cJ^m (X))$. Since
$\epsilon$ is an isomorphism, $T \cong T'$ by the snake lemma. It
follows that $\E_W$ is - as an extension - the pullback of
$\P_W$ via the Abel-Jacobi map and therefore is split if $W$ is
Abel-Jacobi equivalent to zero. This proves the lemma. \hfill
$\square$ \\
As a corollary we get: \\

\proclaim Theorem 1.
$\E_W$ is the pullback of $\P_W$ as extensions via the Abel-Jacobi map.
In particular $\E$ is the pullback of $\P$ as torseurs
via the Abel-Jacobi map.
\par
\ \\ \\ \\
\underline{{\bf \S \ 5. Abel-Jacobi versus Incidence
Equivalence}} \\
\ \\ \\
The notion of ``incidence equivalence'' was introduced by
Griffiths ([Griffiths 69]). Let $X$ be projective, smooth of
dimension $d$ over $\C$. \\
\ \\
{\bf Definition:} A cycle $W \in CH^{i}_{{\rm hom}} (X)$ is called
{\it incidence
equivalent to zero}, i.e. $W \sim_{{\rm inc}} \ 0$, if for all pairs $(T,B)$
with $T$ a smooth, projective variety and $B \in
CH^{d+1 -i} (T \times X)$,
$$
B(W) := (pr_T )_* (B \cap pr^{*}_{X} W) = \{ t \in T | B_t \cap W
\neq \emptyset \}
$$
is linearly equivalent to zero in $\Pic^0 (T)$. \\
\ \\
{\bf Remark:} $B(W)$ is a divisor by a dimension count and in $\Pic^0 (T)$,
since $W \sim_{{\rm hom}} \ 0$. If $i=1$ we get back the notion
of rational equivalence for divisors (= Abel-Jacobi
equivalence). \\
Bloch defines $W \sim_{{\rm inc}} \ 0$ if $\E_W$ is a split
extension. However it is not a priori clear that this gives the
same definition. See the discussion below. \\

\proclaim Lemma 3.
If $W \sim_{A\cJ} \ 0$, then also $W \sim_{{\rm inc}} \ 0$.
\par
\ \\
{\bf Proof.} $B$ is a correspondence
$$
CH^{i}_{{\rm hom}} (X) \to CH^{1}_{{\rm hom}} (T) = \Pic^0 (T)
$$
via $W \mapsto B(W)$. Consider the diagram ($\psi$ = Abel-Jacobi
homomorphism)
$$
\matrix{
CH^{i}_{{\rm hom}} (X) &{\buildrel B \over \to}&
CH^{1}_{{\rm hom}} (T) = \Pic^0(T)\cr
\downarrow  \psi& & \downarrow S \cr
\cJ^i (X)& {\buildrel \cJ B \over \to}&  \cJ^1 (T) }
$$
The claim follows. \hfill $\square$ \\
\ \\
{\bf Remark:} Griffiths has conjectured, that the opposite is true modulo
torsion for cycles algebraicly equivalent to zero:
If $W \sim_{{\rm inc}} \ 0$, then for some $N \in \N$:
$N \cdot W \sim_{A \cJ} \ 0$. This is true for $i=1,d$ without
torsion and for $i=2$ by [Murre 85]. In any case it follows from
the generalized Hodge conjecture of Grothendieck, as we will see
below.  \\
\ \\
In the remaining discussion let $i=d+1-m$ and $W$ a cycle of dimension $m-1$,
hence codimension $d+1-m$ which is homologous to zero.\\
We want to investigate the relation between $W \sim_{{\rm
inc}} \ 0$ and the splitting of $\E_W$ in the sense of [Bloch89].
To do this consider the
action of the correspondence $B$ on the higher Chow groups.
$B(W) \in \Pic^0 (T)$ gives rise to the exact sequence
$$
A_{B(W)} = \frac{CH^t (|B(W)|,1)}{{\rm Im} \  CH^t (T,1)} \to CH^t
(T,|B(W)|,0) \to CH^t (T) \to 0
$$
where $t = \dim T$ and the corresponding extension $\E_{B(W)}$
is sitting in
$$
1 \to \C^* \to \E_{B(W)} \to CH^{t}_{{\rm hom}} (T) \to 0
$$
and is pulled back via the Albanese map from $\P_W$:
$$
1 \to \C^* \to \P_{B(W)} \to {\rm Alb} (T) \to 0.
$$
The correspondence $B$ induces maps $B^{\sharp} = B^{-1} :
CH^{t}_{{\rm hom}} (T) \to CH^{m}_{{\rm hom}} (X)$
then giving rise to a commutative diagram
$$
\matrix{
0 & \to & A_{B(W)} &\to &  CH^t (T, |B(W)|, 0) &\to &  CH^t (T) &\to& 0 \cr
&& \downarrow B^{\sharp} && \downarrow B^{\sharp}&&\downarrow B^{\sharp}
&& \cr
0 &\to& A_W &\to& CH^m (X,| W|, 0) &\to & CH^m (X) &\to & 0}
$$
This can be seen as follows: By Chow's moving lemma, cycles in
$CH^m(X)$ and $CH^m(X,|W|)$ may be assumed to have support disjoint from $W$.
The same holds for $T$ and $B(W)$. Now if a zero cycle on $T$ does not meet
$B(W)$ then $B^\sharp(Z)$ will not meet $W$. This explains the maps in the
diagram. Everything commutes by the functorial properties of \S \ 2. \\
We thus obtain a diagram
$$
\matrix{
1 & \to & \C^* &\to & \E_{B(W)} &\to & CH^{t}_{{\rm hom}} (T) & \to & 0 \cr
&& || && \downarrow B^{\sharp}&& \downarrow  B^{\sharp} &&\cr
1 &\to& \C^* &\to& \E_W &\to& CH^{m}_{{\rm hom}} (X) &\to& 0}
$$
inducing the identity on $\C^*$, where we assumed that $B$ is reduced without
loss
of generality. \\
This is the main input in order to prove: \\

\proclaim Lemma 4.
(a) If $\E_W$ splits then also $\E_{B(W)}$ splits for all pairs $(T,B)$. \\
(b) Assume additionally that $W \sim_{\rm alg} 0$. Then:\\
$\E_{B(W)}$ splits for all pairs $(T,B)$ $\Longleftrightarrow \E_W^{\rm alg}$
splits, where $\E_W^{\rm alg}$ is the subextension
$$
1 \to \C^* \to \E_W^{\rm alg} \to CH^{m}_{{\rm alg}} (X) \to 0
$$
\par
\ \\
{\bf Proof.}
(a) follows directly from the commutative diagram for every pair $(T,B)$.
To prove (b) suppose $W \sim_{\rm alg} 0$ and that $\E_{B(W)}$ splits
for all pairs $(T,B)$. We have to show that $\E_W^{\rm alg}$ splits.
Take $(T,B)$ such that
$$B^{\sharp}: {\rm Alb(T)} \to \cJ^m(X) $$
parametrizes the whole image $\cJ^m_{\rm alg}(X)$
of the Abel-Jacobi map for cycles algebraically
equivalent to zero. $T$ can be chosen to be an abelian variety itself such
that $B^\sharp$ becomes an isogeny. It follows that the subextension
$$
1 \to \C^\ast \to \P^{\rm alg}_W \to \cJ^m_{\rm alg}(X) \to 0
$$
of $\P$ splits by the Deligne cohomology version of the
commutative diagram above.
Now by theorem 1 and its proof, the extension class of
$\E_W^{\rm alg}$ is pulled back from
${\rm Ext}^1(\psi(CH^m_{\rm alg}(X)),\C^*)=
{\rm Ext}^1(\cJ^m_{\rm alg}(X),\C^*)$.
(b) follows.
\hfill $\square$ \\
\ \\
As a corollary we get: \\

\proclaim Theorem 2.
Assume $W \sim_{\rm alg} 0$. Then
$W \sim_{{\rm inc}} \ 0$ if and only if $\E_W^{\rm alg}$ splits.
\par
\ \\
{\bf Proof.} $W \sim_{{\rm inc}} \ 0$ iff $B(W)$ is
linearly equivalent to zero
by definition, hence iff $\P_{B(W)}$ splits. But for divisors the
splittings of $\E$ and $\P$ are equivalent (as one can see by
restricting to suitable curve) and hence this holds if and
only if $\E_{B(W)}$ splits. Hence by lemma 4 (b) we are finished.
\hfill $\square$ \\
\ \\
It remains to discuss {\it Griffiths' conjecture}: Let $W
\sim_{{\rm alg}} \ 0$. If $W \sim_{{\rm inc}} \ 0$ then some
multiple of $W$ is Abel-Jacobi equivalent to zero. With the help
of theorem 2 it is sufficient to show: If $\E_W$ splits, then
for some $N \in \N$ : $\P_{N\cdot W} $ splits. \\

\proclaim Lemma 5.
Suppose the generalized Hodge conjecture (GHC) holds for $\cJ^m
(X)$, i.e. the largest abelian subvariety $\cJ^{m}_{\rm ab} (X)$ of
$\cJ^m(X)$ is
parametrized by algebraic cycles. Then for $W \in CH_{m-1} (X)$,
$W \sim_{{\rm alg}} \ 0$, we have: If $\E_W$ splits, then for
some $N \in \N$, $\P_{N\cdot W}$ splits.
\par
\ \\
{\bf Proof.}
Remark that $\cJ^{m}_{ab} (X)$ and $\cJ^{d+1-m}_{ab} (X)$ are always dual
to each other modulo isogeny, since the dual of a Hodge substructure is
again one. Therefore GHC implies that
$\cJ^{m}_{\rm alg} (X) = \cJ^{m}_{\rm ab} (X)$ and that
$\psi (N\cdot W)$ lies in the dual abelian
variety of $\cJ^{m}_{{\rm alg}} (X) = \psi (CH^{m}_{{\rm alg}}
(X))$ for some $N \in \N$. Hence it is sufficient to show that the extension
$$
1 \to \C^* \to \P_{N \cdot W, {\rm alg}} \to
\cJ^{m}_{{\rm alg}} (X) \to 0
$$
splits. This in turn is implied by the splitting of $\E_W$ and hence of
$\E_{N \cdot W}^{\rm alg}$ by the proof of lemma 4. \hfill $\square$ \\
\ \\
Actually we have proved more: \\

\proclaim Theorem 3.
Let $W \in CH_{m-1} (X) $ be algebraically equivalent to zero.
Then if $\psi (N \cdot W)$ is contained in the dual of
$\cJ^{m}_{{\rm alg}} (X)$,
the splitting of $\E_W$ implies the
splitting of $\P_{N \cdot W}$.
\par
\ \\
A little generalization of the argument in Lemma 5.2. of [Murre85] gives
another proof of Murre's theorem:
\\
\proclaim Corollary(Murre's theorem). \\
Griffiths' conjecture holds for codimension two cycles on a smooth projective
manifold $X$ over $\C$, i.e. for every cycle $W$ on $X$, algebraicly and
incidence equivalent to zero, $\psi(W)$ is torsion in $\cJ^2(X)$.
\par
\ \\
{\bf Proof.}
We have to verify the assumptions of theorem 3. Let $T$ be the
universal cover (tangent space) of $\cJ^{2}_{{\rm alg}} (X)$
and $H^{3}_a(X) = T \oplus \bar T$. This is a
sub Hodge structure of $H^{3}(X,\C)$ contained in
$H^{1,2}(X) \oplus H^{2,1}(X)$.
Now $\psi (W)$ is contained in the dual of $\cJ^{d-1}_{{\rm alg}} (X)$
($d=\dim(X)$) if we can show that the cup product
$$ H^{3}_a(X) \otimes H^{2d-3}_a(X) \to \C $$
has no left kernel. Let $L$ be the Lefschetz operator, i.e. cup product with
the hyperplane class. It is algebraic and hence maps $H^{k}_a(X)$ to
$H^{k+2}_a(X)$. \\
Let us repeat the proof of Murre's Lemma 5.2. for convenience.
Assume there is a nonzero element $v \in  H^{3}_a(X)$
in the left kernel and decompose it as
$v=v_0+Lv_1$ according to Lefschetz decomposition. The
$v_k$ are primitive. In particular $v_1 \in  H^{1}_a(X)$ which is
the whole of $H^{1}(X,\C)$. Hence $Lv_1$ and $v_0$ are in $H^{3}_a(X)$.
Now decompose $v_0=w_1+w_2$ into types $(1,2)$ and $(2,1)$.
Since $ H^{2p-1}_a(X)$ is a sub Hodge structure, $w_1$ and $w_2$
and their complex conjugates are again in $ H^{2p-1}_a(X)$.
By symmetry we may assume $w_1 \neq 0$ unless $v_0 =0$. In the first case set
$z:= L^{d-3}(\bar w_1) \in H^{2d-3}_a(X)$, otherwise decompose $v_1=y_1+y_2$
into types with $y_1 \neq 0$ (by symmetry) and set $z :=L^{d-2}(\bar y_1)$.
Then for type reasons in the first case
$v \cup z= v_0 \cup z =L^{d-3} (w_1 \cup \bar w_1) \neq 0$
by primitivity of $w_1$, a contradiction. In the second case
$v \cup z= Lv_1 \cup z =L^{d-1} (y_1 \cup \bar y_1) \neq 0$
by primitivity of $y_1$, a contradiction. \hfill $\square$ \\
\ \\ \\
{\bf Remark:} It is easy to see that the proof of Griffiths' conjecture
would go through in any codimension $p$, if in the Lefschetz decomposition
of $v \in H^{2p-1}_a(X)$ every term were again in $H^{2p-2k-1}_a(X)$. This
can be verified
for several types of varieties e.g. rational-like varieties and abelian
varieties and complete intersections in these. For a special case see
[Griffiths-Schmid75].\\
\newpage
\noindent
\begin{center}
\underline{{\bf References}}
\end{center}
\vspace{1cm}
\par \noindent
[Bloch 86] \ \ S. Bloch: Algebraic cycles and higher $K$-theory. Adv.
in Math. 61, 167 - 304 (1986) \\
\par \noindent
[Bloch 89] \ \ S. Bloch: Cycles and Biextensions, in Cont. Math. 83,
19 - 30 (1989) \\
\par \noindent
[Carlson 85] \ \ J. Carlson: The geometry of the extension class of a
mixed Hodge structure, in Proc. Bowdoin conference of the AMS Vol.2, 199-222
(1987) \\
\par \noindent
[Esnault-Viehweg 88] \ \ H. Esnault, E. Viehweg: Deligne-Beilinson
Cohomology, Perspectives in Mathematics, Vol. 4, Acad. Press, 43-91 (1988) \\
\par \noindent
[Griffiths 69] \ \ Ph. Griffiths:
Some results on algebraic cycles on algebraic manifolds,
Proc. Bombay conference, 93 - 191 (1969) \\
\par \noindent
[Griffiths-Schmid75] \ \ Ph. Griffiths, W. Schmid: Recent developments in
Hodge theory, in Discrete subgroups of Liegroups, Bombay Coll. 1975, Oxford
Univ. Press, 31-127 (1975) \\
\par \noindent
[Hain 90] \ \ R. Hain: Biextensions and heights associated to curves of
odd genus; Duke Math. J. 61, 859 - 898 (1990) \\
\par \noindent
[Mumford 69] \ \ D. Mumford: Biextensions of formal groups,
Proc. Bombay Colloquium 1968, Oxford Univ. Press, 307-322 (1969) \\
\par \noindent
[Murre] unpublished notes \\
\par \noindent
[Murre 85] \ \ J. Murre: Abel Jacobi equivalence versus incidence equivalence
for algebraic cycles of codimension two, Topology 24 (3), 361 - 367 (1985)
\end{document}